\documentclass[reqno]{amsart}

\usepackage{amssymb,amscd,upref,xspace}

\newtheorem{theorem}{Theorem}[section]
\newtheorem{lemma}[theorem]{Lemma}
\newtheorem{corollary}[theorem]{Corollary}
\newtheorem{proposition}[theorem]{Proposition}

\theoremstyle{definition}
\newtheorem{definition}[theorem]{Definition}
\newtheorem{example}[theorem]{Example}

\theoremstyle{remark}
\newtheorem{remark}[theorem]{Remark}

\numberwithin{equation}{section}
\numberwithin{figure}{section}

\newcommand{\cprime}{\/{\mathsurround=0pt$'$}}

\newcommand{\abs}[1]{\lvert#1\rvert}
\newcommand{\norm}[1]{\lVert#1\rVert}
\newcommand{\dd}[2]{\dfrac{\partial#1}{\partial#2}}
\newcommand{\ddn}[2]{\partial #1/\partial #2}
\newcommand{\R}{\mathcal{R}}
\newcommand{\E}{\mathcal{E}}
\newcommand{\Ei}{\mathcal{E}^{\infty}}
\newcommand{\CL}[1]{\mathcal{C}^{#1}\Lambda}
\newcommand{\xra}{\xrightarrow}
\newcommand{\xla}{\xleftarrow}
\newcommand{\hj}{\bar{\jmath}}
\newcommand{\ji}{\bar{\jmath}_{\infty}}
\newcommand{\hd}{\bar{d}}
\newcommand{\id}{\mathrm{id}}
\newcommand{\J}{\bar{\mathcal{J}}}
\newcommand{\Ji}{{\bar{\mathcal{J}}}^{\infty}}
\newcommand{\F}{\mathcal{F}}
\newcommand{\hL}{\bar{\Lambda}}
\newcommand{\gl}{\mathfrak{gl}}
\newcommand{\sltwo}{\mathfrak{sl}_2}
\newcommand{\hH}{\bar{H}}
\newcommand{\hS}{\bar{\mathrm{S}}}
\newcommand{\CDiff}{\mathcal{C}\mathrm{Dif{}f}}
\newcommand{\Der}{\mathrm{D}}
\newcommand{\CDer}{\mathcal{C}\mathrm{D}}
\newcommand{\Dv}{\mathrm{D^v}}
\newcommand{\Hom}{\mathrm{Hom}}
\newcommand{\hatotimes}{\mathbin{\widehat{\otimes}}}

\DeclareMathOperator{\im}{im}
\DeclareMathOperator{\coker}{coker}

\DeclareMathOperator{\ord}{ord}
\DeclareMathOperator{\ad}{ad}
\DeclareMathOperator{\rank}{rank}

\newcommand{\cm}{$\mathcal{C}$-mod\-ule\xspace}
\newcommand{\cms}{$\mathcal{C}$-mod\-ules\xspace}
\newcommand{\fm}{$\mathcal{F}$-mod\-ule\xspace}
\newcommand{\fms}{$\mathcal{F}$-mod\-ules\xspace}
\newcommand{\cd}{$\mathcal{C}$-dif\-fer\-en\-tial\xspace}
\newcommand{\css}{the Vinogradov $\mathcal{C}$-spec\-tral
                  sequence\xspace}
\newcommand{\kcc}{Krasil\cprime\-shchik's
                  $\mathcal{C}$-co\-ho\-mol\-ogy\xspace}

\begin{document}

\hfill Preprint DIPS 3/98

\hfill math.DG/9803115

\vspace{5cm}

\title{Notes on the horizontal cohomology}

\author{Alexander Verbovetsky}

\address{Moscow State Technical University and The Diffiety Institute,
         Moscow, Russia. Correspondence to: A.~M.~Verbovetsky,
         Profsoyuznaya 98-9-132, 117485 Moscow, Russia}

\email{verbovet@mail.ecfor.rssi.ru}

\thanks{Talk presented at the conference ``Secondary Calculus and
Cohomological Physics'' in Moscow, August 1997.  It will be published
in the Proceedings of this conference (M.~Henneaux, I.~S.
Krasil{\cprime}shchik, and A.~M. Vinogradov, eds.) in Contemporary
Mathematics, Amer. Math. Soc.}

\thanks{This work was supported in part by RFBR grant 97-01-00462
        and INTAS grant 96-0793.}


\begin{abstract}
This paper is devoted to the horizontal (``characteristic'') cohomology
of systems of differential equations. Recent results on computing the
horizontal cohomology via the compatibility complex are generalized.
New results on \css and \kcc are obtained. As an application of general
theory, the examples of an evolution equation and a $p$-form gauge
theory are explicitly worked out.
\end{abstract}

\maketitle

\section*{Introduction}
\label{sec:intr}

Consider a system of differential equations
\begin{equation*}
F_s(x_i,u^j,\dots,u^j_{\sigma},\dotsc)=0,
\tag{$*$}
\end{equation*}
with independent variables $x_i$, unknown functions $u^j$, and
$u^j_{\sigma}= \ddn{^r u^j}{x_{i_1}\dotsm\partial x_{i_r}}$ being their
derivatives, $\sigma= i_1\dotsc i_r$. Let $\F$ be the algebra of
functions in the variables $x_i$, $u^j$, and $u^j_{\sigma}$. Two
functions are said to be equivalent if they coincide when equations
\thetag{$*$} hold. Denote the quotient algebra by $\F_F$. The
horizontal de~Rham complex for system \thetag{$*$} is the exterior
algebra generated by elements of $\F_F$ and differentials $dx_i$, with
the differential being the total exterior differential $\hd(f\,dx_{i_1}
\wedge\dots\wedge dx_{i_k})=\sum_i(df/dx_i)\,dx_i\wedge dx_{i_1}\wedge
\dots\wedge dx_{i_k}$. The cohomology of this complex, called the
\emph{horizontal} (``characteristic'') cohomology, has been studied
and used extensively in the literature (see the bibliography). This
cohomology plays a central role in the Lagrangian formalism (because
the horizontal cohomology group in degree $n$ is the space of actions
of variational problems constrained by equations \thetag{$*$}). The
horizontal cohomology group in degree $n-1$ consists of conservation
laws for the system under consideration. This fact is basic to finding
conservation laws via \css. The horizontal cohomology in degrees less
than $n-1$ has attracted recently a great deal of interest in the
context of gauge theories (for equations without gauge symmetries this
cohomology is trivial by the famous ``two-line theorem'' of Vinogradov
\cite{Vin2,Vin1}), where it appears as a means of calculating the BRST
cohomology \cite{BarnBrandtHenn1}.

For computing the horizontal cohomology there is a general method
based on \css. It can be outlined as follows. The horizontal cohomology
is the term $E_1^{0,\bullet}$ of \css and thereby related to the terms
$E_1^{p,\bullet}$, $p>0$. For each $p$, such a term is also a
horizontal cohomology but with some nontrivial coefficients. The
crucial observation is that the corresponding module of coefficients is
supplied with a filtration such that the differential of the associated
graded complex is linear over the functions. Hence, the cohomology can
be computed algebraically. For $p=1$, this has been done in
\cite{Tsuji1} and the main result is: The term $E_1^{1,\bullet}$
coincides with the homology of the complex dual to the compatibility
complex for the linearization operator of the system under
consideration. In this paper, this result is generalized to $p>1$.  We
compute also \kcc (see \cite{Kras1,Kras2,KrasKers1,KrasKers3,KrasKers2,
KrasKers4}). In particular, the ``$k$-line theorem'' is proved for both
cases. It should be also stressed that our techniques are capable not
only of obtaining vanishing results, but facilitate the computation of
nonzero part as well.

The author's thinking about computing the horizontal cohomology was
especially stimulated by works of I.~S.~Krasil{\cprime}shchik
\cite{Kras1} and M.~Marvan \cite{Marvan1} where the horizontal
cohomology with coefficients in the Cartan forms and zero-curvature
representations was calculated.  By trying to understand their
calculations, the author arrived at the construction of the present
paper.

The ability to calculate the horizontal cohomology is not unique to the
method described above. There is also an approach based on the
Koszul\,-\,Tate resolution. The discussion of this approach does not
enter into the scope of the present paper. It is desired here merely to
point out that both methods take as a starting point the compatibility
complex and furnish equivalent results although the precise
relationships between them have yet to be uncovered.

The paper is organized as follows. We begin with a short introduction
to the \cd calculus on differential equations. In the opening
subsection of Section \ref{sec:hc}, we discuss the notion of a \cm. A
\cm is, by definition, a module equipped with an action of \cd
operators. Thus, \cms serve as modules of coefficients for horizontal
de~Rham complexes. Such complexes are dealt with in Section
\ref{subsec:hdrc}. In Section \ref{subsec:cc}, we relate the horizontal
cohomology to the cohomology of the compatibility complex. The results
of this section play a key role in the explicit calculation of the
horizontal cohomology. Finally, in Section \ref{sec:appl} we turn to
the examples of computing \css and \kcc.

The author would like to thank Dmitri Gessler, Joseph~S.
Krasil\cprime\-shchik, and Alexander~M. Vinogradov for many useful
discussions and collaborations. It is a pleasure to thank the
organizers of the International Conference on Secondary Calculus and
Cohomological Physics in Moscow for the invitation, providing a
stimulating environment, and financial help. He is also grateful to the
Scuola Internazionale Superiore di Studi Avanzati in Trieste for its
kind hospitality during the author's visit in 1993--95 when a part of
this research was done.

\section{$\protect\mathcal{C}$-differential calculus}\label{sec:cdc}

\subsection{Jet spaces}\label{subsec:jets}
In this section a brief discussion of the \cd calculus on differential
equations is given, further details being available in
\cite{Symm,KLV,Vin5,Olver}.

Let $\pi\colon E\to M$ be an $n+m$-dimensional vector bundle over an
$n$-dimensional base manifold $M$, and let $\pi_\infty\colon J^\infty
(\pi)\to M$ be the infinite jet bundle of local sections of $\pi$.
Denote by $j_\infty(s)(x)$ the infinite jet of a local section $s$ of
the bundle $\pi$ at point $x\in M$. Thus each section $s\colon M\to E$
gives rise to the section $j_\infty(s)\colon M\to J^\infty(\pi)$.

A coordinate system $(x_i,u^j)$, $i=1,\dots,n$, $j=1,\dots,m$, of
$\pi\colon E\to M$ induces coordinates
$(x_i,u^j,u^j_i,\dots,u^j_{\sigma},\dotsc)$ on $J^\infty(\pi)$ such
that
\[u^j_{\sigma}(j_\infty(s)(x))=\dd{^{\abs{\sigma}}s^j}{x_{i_1}\partial
x_{i_2}\dotsm\partial x_{i_r}},\]
where $\sigma=i_1 i_2 \dotsc i_r$ is a multi-index, $\abs{\sigma}=r$.

Denote the algebra of functions on $J^\infty(\pi)$ by $\F=\F(\pi)$.
Functions on $J^\infty(\pi)$ are smooth functions of a \emph{finite}
number of coordinates. Let $\alpha\colon W\to J^\infty(\pi)$ be a
vector bundle and $\Gamma(\alpha)$ be the set of local sections of
$\alpha$.  Naturally, $\Gamma(\alpha)$ is a module over the algebra
$\F$. In the sequel, by a module we always mean a module of this form.
Recall that homomorphisms of modules are in one-to-one correspondence
with morphisms of bundles, so we can and will not distinguish between
bundles and modules.

On an infinite jet space, there exists the \emph{\cd calculus}. In
coordinate language, this is the total derivatives calculus. Let us
describe \cd operators, corresponding jets, and forms in more details.

Consider two \fms $P$ and $P_1$. A differential operator $\Delta\colon
P\to P_1$ is called \emph{$\mathcal{C}$-dif\-fer\-en\-tial}, if it can
be restricted to the manifolds of the form $j_\infty(s)(M)$, where $s$
is a local section of $\pi$. In other words, $\Delta$ is a \cd
operator, if the equality $j_\infty(s)^*(\varphi)=0$, $\varphi\in P$,
implies $j_\infty(s)^*(\Delta(\varphi))=0$.

In local coordinates, \cd operators have the form
\[\begin{pmatrix}
\sum_{\sigma}a^\sigma_{11}D_\sigma & \dots &
\sum_{\sigma}a^\sigma_{1m_1}D_\sigma \\
\hdotsfor{3}\\
\sum_{\sigma}a^\sigma_{m_2 1}D_\sigma & \dots &
\sum_{\sigma}a^\sigma_{m_2 m_1}D_\sigma
\end{pmatrix},\]
where $a^{\sigma}_{ij}\in\F$, $D_{\sigma}=D_{i_1}\circ\dots\circ
D_{i_r}$ for $\sigma=i_1\dotsc i_r$, and
$D_i=\ddn{}{x_i}+\sum_{j,\sigma}u_{\sigma i}^j \ddn{}{u_{\sigma}^j}$ is
the $i$-th \emph{total derivative operator}.

We shall denote the \fm of \cd operators from $P$ to $P_1$ by
$\CDiff(P,P_1)$. In this module, there exists a filtration by the
modules $\CDiff_k(P,P_1)$ that consist of \cd operators of order $\le
k$.

Next, define the \emph{module of horizontal jets}. Let $P$ be an \fm.
Say that two elements of $P$ are horizontally equivalent up to order $k
\le\infty$ at point $\theta\in J^\infty(\pi)$, if their total
derivatives up to order $k$ coincide at $\theta$. The horizontal jet
space $\bar{J}_\theta^k (P)$ is $P$ modulo this relation, and the
collection $\bar{J}^k(P)=\bigcup
_{\theta\,\in\,J^\infty(\pi)}\bar{J}^k_\theta(P)$ constitutes the
\emph{horizontal jet bundle} $\bar{J}^k(P)\to J^\infty(\pi)$. We denote
the module of sections of horizontal jet bundle by $\J^k(P)$.

As with the usual jet bundles, there exist the natural \cd operators
\[\hj_k\colon P\to\J^k(P),\]
 and the natural projections $\nu_{k,l}\colon\J^k(P)\to\J^l(P)$ such
that $\nu_{k,l}\circ\hj_k=\hj_l$. For any \cd operator $\Delta\colon
P\to P_1$ of order $\le k$, there is a unique $\F$-homomorphism
$\varphi_\Delta\colon\J^k(P)\to P_1$ such that
$\Delta=\varphi_\Delta\circ\hj_k$. The composition
\[\hj_l\circ\Delta\in\CDiff_{k+l}(P,\J^l(P_1))\]
is called the $l$-th prolongation of $\Delta$ and the corresponding
homomorphism from $\J^{k+l}(P)$ to $\J^l(P_1)$ is denoted by $\varphi_\Delta
^l$. In particular, $\Delta$ gives rise to the mapping $\varphi_\Delta^\infty
\colon\Ji(P)\to\Ji(P_1)$.

Let $\Der(\pi)$ be the \fm of vector fields on $J^\infty(\pi)$. Define
the \emph{module of Cartan vector fields} as the intersection
$\CDer(\pi)=\Der( \pi)\cap\CDiff_1(\F,\F)$. In coordinates, a Cartan
vector field is $\sum_i f_i D_i$, $f_i\in\F$. A vector field
$X\in\Der(\pi)$ is called \emph{vertical} if $X(j_\infty(s)^*(h))=0$
for all functions $h\in C^\infty (M)$. Locally a vertical vector field
has the form $\sum_{j,\sigma} f_{j,\sigma}\ddn{}{u_\sigma^j}$. Denote
the module of vertical vector fields by $\Dv(\pi)$. Clearly,
$\Der(\pi)=\CDer(\pi)\oplus\Dv(\pi)$.

Further, consider the module $\Lambda^k(\pi)$ of differential $k$-forms on
$J^\infty(\pi)$. A form $\omega\in\Lambda^k(\pi)$ is called a \emph{Cartan
form}, if $\omega$ satisfies $j_\infty(s)^*(\omega)=0$ for every local
section $s$ of $\pi$. Equivalently, Cartan forms are forms that vanish
on the Cartan vector fields. The set of all Cartan forms defines an
ideal \[\CL{}^*(\pi)=\bigoplus_{k\,\ge\,0}\CL{}^k(\pi)\] in the ring
\[\Lambda^*(\pi)=\bigoplus_{k\,\ge\,0}\Lambda^k(\pi)\]
of all forms on $J^\infty(\pi)$. In coordinates, this ideal is generated by
the Cartan $1$-forms
\[\omega^j_{\sigma}=du^j_{\sigma}-\sum_i u^j_{\sigma i}\,dx_i\]
for all $j$ and $\sigma$. The quotient algebra
\[\hL^*(\pi)=\Lambda^*(\pi)/\CL{}^*(\pi)\]
is called the algebra of \emph{horizontal forms} on $J^\infty(\pi)$. In
coordinates, a horizontal $k$-form is a sum of terms of the form
\[f\,dx_{i_1}\wedge\dots\wedge dx_{i_k},\quad f\in\F(\pi).\]
The exterior derivative
\[d\colon\Lambda^k(\pi)\to\Lambda^{k+1}(\pi)\]
gives rise to the \emph{horizontal differential}
\[\hd\colon\hL^k(\pi)\to\hL^{k+1}(\pi),\]
since the ideal of Cartan forms $\CL{}^*(M)$ is stable with respect to $d$:
$d(\CL{}^*(\pi))\subset\CL{}^*(\pi)$. Thus we get the \emph{horizontal
de~Rham complex}
\[0 \xra{} \F(\pi) \xra{\hd} \hL^1(\pi) \xra{\hd}
\hL^2(\pi) \xra{\hd} \dotsb \xra{\hd} \hL^n(\pi) \xra{} 0.\]

More general, consider the filtration in the de~Rham complex on the jet space
$J^\infty(\pi)$
\[\dotsb \subset \CL{k+1}^* \subset \CL{k}^* \subset \dotsb \subset
\CL{2}^* \subset \CL{}^* \subset \Lambda^*,\]
where $\CL{k}^*$ is the $k$-th exterior power of the Cartan ideal $\CL{}^*$.
Evidently, $d(\CL{k}^*)\subset\CL{k}^*$, so that we obtain a spectral
sequence converging to the de~Rham cohomology of $J^\infty(\pi)$. The zero
term of the spectral sequence is $E_0^{p,q}=\CL{p}^{p+q}=\CL{}^1\wedge\dotsm
\wedge\CL{}^1\otimes\hL^q$. This spectral sequence is called \css \cite{Vin3,
Vin2}. Further details on this important spectral sequence can be found in
\cite{Vin1,Symm,Vin4,Tsuji2,Tsuji3,Tsuji1,BryantGriff1,Ander1,Ander2,%
VerbVinGess,Marvan4} and below in this paper as well.

Given a \cd operator $\Delta\in\CDiff_k(P,P_1)$, we define the ($l$-th)
symbol $\sigma(\Delta)$ of $\Delta$ by the following commutative diagram
\[\begin{CD}
0 @>>> S^{k+l}(\hL^1)\otimes P @>>> \J^{k+l}(P) @>\nu_{k+l,k+l-1}>>
\J^{k+l-1}(P) @>>> 0 \\
@. @VV\sigma(\Delta)V @VV\varphi_{\Delta}^{k+l}V
@VV\varphi_{\Delta}^{k+l-1}V @. \\
0 @>>> S^k(\hL^1)\otimes P_1 @>>> \J^k(P_1) @>\nu_{k,k-1}>>
\J^{k-1}(P_1) @>>> 0.
\end{CD}\]
The rows of the diagram are exact, with inclusions $S^k(\hL^1)\otimes P\to
\J^k(P)$ given by $df_1\dotsm df_k\otimes p\mapsto [\dotsc[\hj_k,f_1],\dots,
f_k](p)$, where $p\in P$, $[\cdot\,,f_i]$ is the commutator with the operator
of multiplication by $f_i\in\F$.

Let us consider the pullback of $\pi$ along the projection $J^\infty(\pi)\to
M$ and denote the module of sections of this vector bundle by $\varkappa$. It
is readily seen that for any point $\theta=j_\infty(s)(x)\in J^\infty(\pi)$
one has
\begin{equation}\label{eq:ci}
T_{\theta}(J^\infty(\pi))=J_x^\infty(\pi)=\bar{J}_{\theta}^\infty(\varkappa).
\end{equation}
This yields the canonical isomorphism
\begin{equation}\label{eq:dvhj}
\Dv(\pi)=\Ji(\varkappa).
\end{equation}
The dual isomorphism reads
\begin{equation}\label{eq:cfcd}
\CL{}^1(\pi)=\CDiff(\varkappa,\F).
\end{equation}
In coordinates, the form $\omega_{\sigma}^j$ under this isomorphism is the
operator $(\dots,D_\sigma,\dotsc)$, with $D_{\sigma}$ on $j$-th place.

It is clear that the Cartan $k$-forms can be identified with multilinear
skew-symmetric \cd operators in $k$ arguments.

\subsection{Differential equations}\label{subsec:de}
Pick up a system of $k$-th order partial differential equations
\begin{equation}\label{eq:fz}
F_s(x_i,u^j,\dots,u^j_{\sigma},\dotsc)=0,\quad s=1,\dots,l.
\end{equation}
We shall consider $F=(F_1,\dots,F_l)$ as an element of an \fm $P$. The system
$F=0$ defines a subbundle $\E\to M$ of the jet bundle $J^k(\pi)\to M$. A
solution to the system of differential equations is a section $s\colon M\to
E$ such that $j_k(s)(M)\subset\E$. One can define the \emph{infinite
prolongation} $\Ei$ of $\E$ by equations $F_s=0$, \ldots, $D_{\sigma}F_s=0$,
\ldots\ This system is equivalent to \eqref{eq:fz} in the sense that both
have the same set of solutions. For brevity, $\E$ and $\Ei$ will be referred
to as ``equation''.

We define the function algebra $\F(\E)$ on $\Ei$ to be the restriction of the
algebra $\F(\pi)$ on $J^\infty(\pi)$. It is straightforward to show that one
can pull all ingredients of the \cd calculus on $J^\infty(\pi)$,
discussed in the previous subsection, back to $\Ei$. Thus, we have only
one thing to do here:  to generalize formulae \eqref{eq:dvhj} and
\eqref{eq:cfcd}.

Pick a point $\theta=j_\infty(s)(x)\in\Ei$. In view of \eqref{eq:ci}, the
tangent space $T_{\theta}(\Ei)$ is isomorphic to a subspace $R_\theta\subset
\bar{J}^\infty(\varkappa)$. In coordinates, the subspace $T_\theta(\Ei)
\subset T_\theta(J^\infty(\pi))$ is given by the equations $\sum_{j,\sigma}
(\ddn{D_{\tau}(F_s)}{u_{\sigma}^j})\,du_\sigma^j=0$, $\tau$ being a
multi-index. Therefore $R_\theta$ is defined by
\begin{equation}\label{eq:sdu}
\sum_{j,\sigma}\dd{D_{\tau}(F_s)}{u_{\sigma}^j}w_{\sigma}^j=0,
\end{equation}
where $w_{\sigma}^j=\pi_\infty^*(u_{\sigma}^j)$ are coordinates on $\bar{J}^
\infty(\varkappa)$. Clearly, system \eqref{eq:sdu} can be rewritten in the
form
\[\sum_{j,\sigma}D_{\tau}(\dd{F_s}{u_{\sigma}^j}w_{\sigma}^j)=0\]
or $D_{\tau}(\ell_F(w^j))=0$, where
\[\ell_F=\begin{pmatrix}
\sum_{\sigma}\dd{F_1}{u_{\sigma}^1}D_{\sigma}& \dots &
\sum_{\sigma}\dd{F_1}{u_{\sigma}^m}D_{\sigma} \\
\hdotsfor[2]{3} \\
\sum_{\sigma}\dd{F_l}{u_{\sigma}^1}D_{\sigma}& \dots &
\sum_{\sigma}\dd{F_l}{u_{\sigma}^m}D_{\sigma}
\end{pmatrix}\]
is the \emph{operator of universal linearization} for $F$.

Now our discussion can be summarized as follows.

\begin{proposition}\label{prop:ker}
For any differential equation $\E$,
\begin{enumerate}
\item The module $\Dv(\E)$ is isomorphic to the kernel of the homomorphism
$\varphi_{\ell_F}^\infty\colon$ $\Ji(\varkappa)\to\Ji(P)$.
\item The module $\CL{}^1(\E)$ is isomorphic to $\CDiff(\varkappa,\F)$
modulo the submodule consisting of the operators of the form
$\nabla\circ\ell_F$, $\nabla\in\CDiff(P,\F)$.
\end{enumerate}
\end{proposition}
We use the notation $\Dv(\E)$, $\CL{k}^s(\E)$, and so on, for the
corresponding modules on $\Ei$.

\begin{remark}
The constructions and results covered in this paper are valid not only for
equations $\Ei$, but for arbitrary diffieties \cite{Vin4,Symm} as well.
Recall, that a diffiety is an infinite-dimensional manifold furnished
with an involutive finite-dimensional distribution that locally is of
the form $\Ei$ endowed with the distribution of Cartan fields.
\end{remark}

\section{Horizontal cohomology}\label{sec:hc}

\subsection{$\protect\mathcal{C}$-modules on differential equations}
\label{subsec:c-mod}

Fix an equation $\Ei$. Let $\F=\F(\E)$ be the algebra of functions on $\Ei$.
An \fm $Q$ is called a \emph{$\mathcal{C}$-mod\-ule}, if $Q$ is endowed
with a left module structure over the ring $\CDiff(\F,\F)$, i.e., for
any scalar \cd operator $\Delta\in\CDiff_k(\F,\F)$ there exists an
operator $\Delta_Q\in \CDiff_k(Q,Q)$ with

\begin{enumerate}
\item $(\sum_i f_i\Delta_i)_Q=\sum_i f_i(\Delta_i)_Q,\quad f_i\in\F$,
\item $(\id_{\F})_Q=\id_Q$,
\item $(\Delta_1\circ\Delta_2)_Q=(\Delta_1)_Q\circ(\Delta_2)_Q$.
\end{enumerate}

In other words, a \cm is a module equipped with a \emph{flat horizontal
connection}, i.e., with an action on $Q$ of the module $\CDer=\CDer(\E)$ of
Cartan vector fields, $X\mapsto\nabla_X$, that is $\F$-linear:
\[\nabla_{fX+gY}=f\nabla_X+g\nabla_Y,\quad f,g\in\F,\quad X,Y\in\CDer,\]
satisfies the Leibnitz rule:
\[\nabla_X(fq)=X(f)q+f\nabla_X(q),\quad
q\in Q,\quad X\in\CDer,\quad f\in\F,\]
and is a Lie algebra homomorphism:
\[[\nabla_X,\nabla_Y]=\nabla_{[X,Y]}.\]

The coordinate description of a flat horizontal connection looks as
\[\nabla_{D_i}(s_j)=\sum_k\Gamma^k_{ij}s_k,\quad\Gamma^k_{ij}\in\F,\]
where $s_j$ are basis elements of $Q$.

\begin{remark}
Let $Q$ be the module of sections of a vector bundle $\tau\colon W\to\Ei$,
$Q=\Gamma(\tau)$. A flat horizontal connection on $Q$ defines a completely
integrable $n$-dimensional linear distribution on $W$ that is projected onto
the Cartan distribution on $\Ei$. Thus, geometrically, a \cm is the module of
sections of a linear covering (see \cite{Symm,KrasVin1}).

In coordinates the covering has the form
\[\widetilde{D_i}=D_i+\sum_{j,k}\Gamma^k_{ij}w^j\dd{}{w_k},\]
where $w^i$ are fiber coordinate on $W$.
\end{remark}

Here are basic examples of \cms.
\begin{example}
The simplest example of a \cm is $Q=\F$ with the usual action of horizontal
operators.
\end{example}

\begin{example}
The module of vertical vector fields $Q=\Dv=\Dv(\E)$ with the connection
\[\nabla_X(Y)=[Y,X],\quad X\in\CDer,\quad Y\in\Dv\]
is a \cm.
\end{example}

\begin{example}
Next example is the modules of Cartan forms
$Q=\CL{k}^k=\CL{k}^k(\E)$. A vector field $X\in\CDer$ acts on $\CL{k}^k$
as the Lie derivative $L_X$. It is easily seen that in coordinates we have
\[(D_i)_{\CL{k}^k}(\omega_{\sigma}^j)=\omega_{\sigma i}^j\]
\end{example}

\begin{example}
The infinite jet module $Q=\Ji(P)$ of an \fm $P$ is a \cm via
\[\nabla_X(f\ji(p))=X(f)\ji(p),\]
where $X\in\CDer$, $f\in\F$, $p\in P$.
\end{example}

\begin{example}
Let us dualize the previous example. It is clear that for any \fm $P$ the
module $Q=\CDiff(P,\F)$ is a \cm. The action of horizontal operators is the
composition.
\end{example}

\begin{example}\label{exmp:ker}
More generally, let $\Delta\colon P\to P_1$ be a \cd operator and $\varphi_
\Delta^\infty\colon \Ji(P)\to\Ji(P_1)$ be the corresponding prolongation of
$\Delta$. Obviously, $\varphi_\Delta^\infty$ is a morphism of \cms, i.e., a
homomorphism over the ring $\CDiff(\F,\F)$, so that $\ker\varphi_\Delta^
\infty$ and $\coker\varphi_\Delta^\infty$ are \cms.

Dually, the operator $\Delta$ gives rise to the morphism of \cms
$\CDiff(P_1,\F)\to\CDiff(P,\F)$, $\nabla\mapsto\nabla\circ\Delta$. Thus the
kernel and cokernel of this map are \cms.
\end{example}

\begin{example}
Given two \cms $Q_1$ and $Q_2$, we can define \cm structures on $Q_1\otimes_
{\F} Q_2$ and $\Hom_{\F}(Q_1,Q_2)$ by
\begin{align*}
\nabla_X(q_1\otimes q_2)&=\nabla_X(q_1)\otimes q_2+q_1\otimes\nabla_X(q_2),\\
\nabla_X(f)(q_1)&=\nabla_X(f(q_1))-f(\nabla_X(q_1)),
\end{align*}
where $X\in\CDer$, $q_1\in Q_1$, $f\in\Hom_{\F}(Q_1,Q_2)$.

For instance, we have \cm structures on $Q=\Ji(P)\otimes_{\F}\CL{k}^k$ and
$Q=\CDiff(P,\CL{k}^k)$ for any \fm $P$.
\end{example}

\begin{example}
Let $\mathfrak{g}$ be a Lie algebra and $\rho\colon\mathfrak{g}\to\gl(W)$ a
linear representation of $\mathfrak{g}$. Each $\mathfrak{g}$-valued
horizontal form $\omega\in\hL^1(\E){\otimes}_{\mathbb{R}}\,\mathfrak{g}$ that
satisfies the horizontal Maurer\,-\,Cartan condition
$\hd\omega+\frac12[\omega,\omega]=0$ defines on the module $Q$ of
sections of the trivial vector bundle $\Ei\times W\to\Ei$ the following
\cm structure:  \[\nabla_X(q)_a=X(q)_a+\rho(\omega(X))(q_a),\] where
$X\in\CDer$, $q\in Q$, $a\in\Ei$. Such \cms are called
\emph{zero-curvature representations} over $\Ei$ (cf.\ \cite{Marvan1}). Take
the example of the KdV equation (in the form $u_t=uu_x+u_{xxx}$) and
$\mathfrak{g}=\sltwo(\mathbb{R})$. Then there exists a one-parameter family
of Maurer\,-\,Cartan forms $\omega(\lambda)=A_1(\lambda)\,\hd x+A_2(\lambda)
\,\hd t$, $\lambda$ being a parameter:
\[A_1(\lambda)=\begin{pmatrix}
0& -(\lambda+u)\\
\dfrac16& 0
\end{pmatrix},
\quad A_2(\lambda)=
\begin{pmatrix}
-\dfrac16 u_x&
-u_{xx}-\dfrac13 u^2+\dfrac13 \lambda u+\dfrac23 \lambda^2\\
\dfrac1{18}u-\dfrac19 \lambda&
\dfrac16 u_x
\end{pmatrix}.\]
This is the zero-curvature representation used in the inverse scattering
method.

In coordinates, if the form $\omega$ is given by $\omega=\sum_iA_i\,dx_i$,
$A_i\in\mathfrak{g}$, then for any \cd operator $\Delta$ the coordinate
description of the operator $\Delta_Q$ can be obtained by replacing all
occurrences of $D_i$ with $D_i+\ad A_i$.
\end{example}

\begin{remark}
In parallel with left \cms one can consider \emph{right} \cms, i.e., right
modules over the ring $\CDiff(\F,\F)$. There is a natural way to pass from
left \cms to right ones and back. Namely, for any left module $Q$ set
\[\mathrm{B}(Q)=Q\otimes_{\F}\hL^n(\E),\]
with the right action of $\CDiff(\F,\F)$ on $\mathrm{B}(Q)$ given by
\begin{align*}
(q\otimes\omega)f&=fq\otimes\omega=q\otimes f\omega,\quad f\in\F,\\
(q\otimes\omega)X&=-\nabla_X(q)\otimes\omega-q\otimes L_X\omega,\quad
X\in\CDer.
\end{align*}
One can easily verify that $\mathrm{B}$ determines an equivalence between the
categories of left \cms and right \cms.
\end{remark}

Take a \cm $Q$. By definition, for a scalar \cd operator $\Delta\colon\F\to
\F$ there exists the operator $\Delta_Q\colon Q\to Q$. In fact one has more:

\begin{proposition}
Let $P,S$ be $\F$-modules. Then there exists a unique mapping
\[\CDiff_k(P,S)\to\CDiff_k(P\otimes_{\F}Q,S\otimes_{\F}Q),\qquad
\Delta\mapsto\Delta_Q,\]
such that the following conditions hold:
\begin{enumerate}
\item if $P=S=\F$ then the mapping is given by the \cm structure on $Q$,
\item $(\sum_i f_i\Delta_i)_Q=\sum_i f_i(\Delta_i)_Q,\quad f_i\in\F$,
\item if $\Delta\in\CDiff_0(P,S)=\Hom_{\F}(P,S)$ then
      $\Delta_Q=\Delta\otimes_{\F}\id_Q$,
\item if $R$ is another $\F$-module and $\Delta_1:P\to S,\Delta_2:S\to R$
are \cd operators, then $(\Delta_2\circ\Delta_1)_Q=(\Delta_2)_Q\circ
(\Delta_1)_Q$.
\end{enumerate}
\end{proposition}
\begin{proof}\label{prop:mq}
The uniqueness is obvious. To prove the existence consider the family of
operators $\Delta(p,s^*)\colon\F\to\F$, $p\in P$, $s^*\in S^*=\Hom_{\F}(S,
\F)$, $\Delta(p,s^*)(f)=s^*(\Delta(fp))$, $f\in\F$. Clearly, the operator
$\Delta$ is defined by the family $\Delta(p,s^*)$. The following statement is
also obvious.
\begin{lemma}[\cite{Vin1}]
For the family of operators $\Delta[p,s^*]\in\CDiff_k(\F,\F)$, $p\in P$, $s^*
\in S^*$, we can find an operator $\Delta\in\CDiff_k(P,S)$ such that $\Delta
[p,s^*]=\Delta(p,s^*)$, if and only if
\begin{align*}
\Delta[p,\sum_i f_is_i^*]&=\sum_i f_i\Delta[p,s_i^*],\\
\Delta[\sum_i f_ip_i,s^*]&=\sum_i\Delta[p_i,s^*]f_i.
\end{align*}
\end{lemma}

In view of this lemma, the family of operators
\[\Delta_Q[p\otimes q,s^*\otimes q^*](f)=q^*(\Delta(p,s^*)_Q(fq))\]
uniquely determines the operator $\Delta_Q$.
\end{proof}

\subsection{The horizontal de~Rham complex}\label{subsec:hdrc}
Consider a complex of \cd operators
$\dotsb\xra{}P_{i-1}\xra{\Delta_i}P_i\xra {\Delta_{i+1}} P_{i+1} \xra{}
\dotsb$. Multiplying it by a \cm $Q$ and taking into account
Proposition \ref{prop:mq}, we obtain the complex
\[\dotsb \xra{} P_{i-1}\otimes Q \xra{(\Delta_i)_Q} P_i\otimes Q
\xra{(\Delta_{i+1})_Q} P_{i+1}\otimes Q \xra{} \dotsb.\]
Applying this construction to the horizontal de~Rham complex, we get
\emph{horizontal de~Rham complex with coefficients in $Q$}:
\[0 \xra{} Q \xra{\hd_Q} \hL^1\otimes_{\F}Q \xra{\hd_Q} \dotsb
\xra{\hd_Q} \hL^n \otimes_{\F}Q \xra{} 0,\]
where $\hL^i=\hL^i(\E)$.
The differential $\hd=\hd_Q$ can also be defined by
\begin{align*}
(\hd q)(X)&=\nabla_X(q),\quad q\in Q,\\
\hd(\omega\otimes q)&=\hd\omega\otimes q+(-1)^p\omega\wedge\hd q,\quad
\omega\in\hL^p.
\end{align*}
One easily sees that a morphism $f\colon Q_1\to Q_2$
of \cms gives rise to a chain mapping of the de~Rham complexes:
\[\begin{CD}
0 @>>> Q_1 @>\hd>> \hL^1\otimes_{\F}Q_1 @>\hd>> \dotsb @>\hd>>
\hL^n\otimes_{\F}Q_1 @>>> 0 \\
@.   @VVV   @VVV   @.   @VVV   @. \\
0 @>>>Q_2 @>\hd>> \hL^1\otimes_{\F}Q_2 @>\hd>> \dotsb @>\hd>>
\hL^n\otimes_{\F}Q_2 @>>> 0.
\end{CD}\]
The cohomology of the horizontal de~Rham complex with coefficients in
$Q$ is said to be \emph{horizontal cohomology} and is denoted by
$\hH^i(Q)$.

Let us discuss some examples of the horizontal de~Rham complexes.

\begin{example}
The horizontal de~Rham complex with coefficients in $\Ji(P)$
\[
0 \xra{} \Ji(P) \xra{\hd} \hL^1\otimes\Ji(P) \xra{\hd}
\hL^2\otimes\Ji(P) \xra{\hd} \dotsb \xra{\hd} \hL^n\otimes\Ji(P)\xra{}0
\]
turns out to be the project limit of the \emph{horizontal Spencer
complexes}
\begin{equation}\label{eq:hsp}
0 \xra{} \J^k(P) \xra{\bar{S}} \hL^1\otimes\J^{k-1}(P) \xra{\bar{S}}
\hL^2\otimes\J^{k-1}(P) \xra{\bar{S}}\dotsb \xra{\bar{S}}
\hL^n\otimes\J^{k-n}(P) \xra{} 0,
\end{equation}
where $\bar{S}(\omega\otimes\hj_l(p))=\hd\omega\otimes\hj_{l-1}(p)$. As usual
Spencer complexes, they are exact in positive degrees and
\[H^0(\hL^{\bullet}\otimes\J^{k-\bullet}(P))=P.\]
The proof is standard but we sketch it now because it involves issues which
will be needed further. Consider the short exact sequence of the complexes in
Diagram \ref{eq:dcs}, with $\hS^k=\mathrm{S}^k(\hL^1)$ being the symmetric
power of $\hL^1$.
\begin{figure}\[\begin{CD}
@. 0 @. 0 @. 0 @.  \\
@. @VVV @VVV @VVV @.\\
0 @>>> \hS^k\otimes P @>>> \J^k(P) @>>> \J^{k-1}(P)
@>>> 0 \\
@.  @VV\bar{\delta} V @VV \bar{S} V @VV \bar{S} V @. \\
0 @>>> \hL^1\otimes\hS^{k-1}\otimes P @>>>
\hL^1\otimes\J^{k-1}(P) @>>> \hL^1\otimes\J^{k-2}(P) @>>> 0\\
@.  @VV\bar{\delta} V @VV \bar{S} V @VV \bar{S} V @. \\
0 @>>> \hL^2\otimes\hS^{k-2}\otimes P @>>>
\hL^2\otimes\J^{k-2}(P) @>>> \hL^2\otimes\J^{k-2}(P) @>>> 0\\
@.  @VV\bar{\delta} V @VV \bar{S} V @VV \bar{S} V @. \\
@. \vdots @. \vdots @. \vdots @.
\end{CD}\]\caption{}\label{eq:dcs}\end{figure}
The left column, called the \emph{horizontal Spencer $\delta$-complex}, is a
complex of homomorphisms. The operator $\bar{\delta}\colon\hL^s\otimes\hS^r
\otimes P\to\hL^{s+1}\otimes\hS^{r-1}\otimes P$ is defined by $\bar{\delta}(
\omega\otimes u\otimes p)=(-1)^s \omega\wedge i(u)\otimes p$, where $i\colon
\hS^r\to\hL^1\otimes\hS^{r-1}$ is the natural inclusion. Dropping the
multiplier $P$ and considering the sequence at a point of $\Ei$, we get the
de~Rham complex with polynomial coefficients. This proves that the Spencer
$\delta$-complex is exact. Hence the second and third columns have the same
cohomology, so that the cohomology of \eqref{eq:hsp} can be computed at $k=0$
and the desired statement is proved.\qed

\begin{remark}
Here we encounter the phenomenon mentioned in the Introduction: the existence
of a filtration in the \cm of coefficients, such that the associated graded
complex is a complex of $\F$-homomorphisms, makes it possible to calculate
the horizontal cohomology \emph{locally}.
\end{remark}

Now, let us multiply the previous diagram by a \cm $Q$ (possibly of infinite
rank). Arguing as above, we see that the complex
\[0 \xra{} \Ji(P)\hatotimes Q \xra{\hd} \hL^1\otimes\Ji(P)
\hatotimes Q \xra{\hd} \dotsb \xra{\hd} \hL^n\otimes\Ji(P)
\hatotimes Q \xra{} 0\]
is exact in positive degrees and
\[H^0(\hL^{\bullet}\otimes\Ji(P)\hatotimes Q)= P\otimes Q.\]
Here
\[\Ji(P)\hatotimes Q=\projlim\J^k(P)\otimes Q.\]
\end{example}

\begin{example}
The dualization of the previous example is the following. The
coefficient module is $\CDiff(P,\F)$. The corresponding horizontal
de~Rham complex multiplied by a \cm $Q$ has the form
\[
0 \xra{} \CDiff(P,\F)\otimes Q \xra{\hd} \CDiff(P,\hL^1)\otimes Q
\xra{\hd} \dotsb \xra{\hd} \CDiff(P,\hL^n)\otimes Q \xra{} 0.
\]
As in the previous example, it is easily shown that
\begin{align*}
H^i(\CDiff(P,\hL^{\bullet})\otimes Q)&=0\quad\text{for $i<n$,} \\
H^n(\CDiff(P,\hL^{\bullet})\otimes Q)&=\hat{P}\otimes Q,
\end{align*}
where $\hat{P}=\Hom_{\F}(P,\hL^n)$.

One can use this fact to define the notion of adjoint operator as follows
(cf.\ \cite{Tsuji1,Verb1}). Any horizontal differential operator
\[\Delta\colon P\to P_1\]
gives rise to the morphism of complexes
\[\begin{CD}
       0             @.         0        \\
      @VVV                    @VVV \\
\CDiff(P_1,\F)    @>>> \CDiff(P,\F) \\
      @VVV                    @VVV \\
\CDiff(P_1,\hL^1) @>>> \CDiff(P,\hL^1) \\
      @VVV                    @VVV \\
\CDiff(P_1,\hL^2) @>>> \CDiff(P,\hL^2) \\
      @VVV                    @VVV \\
      \vdots        @.        \vdots \\
      @VVV                    @VVV \\
\CDiff(P_1,\hL^n) @>>> \CDiff(P,\hL^n) \\
      @VVV                    @VVV \\
       0             @.         0
\end{CD}\]
and hence the map of the cohomology groups $\Delta^*\colon\hat{P}_1\to
\hat{P}$, which is the \emph{adjoint operator}. The reader will have no
difficulty in showing that
\begin{enumerate}
\item $(\Delta_1\circ\Delta_2)^*=\Delta_2^*\circ\Delta_1^*$;
\item if $\Delta\in\CDiff_k(P,P_1)$ then
$\Delta^*\in\CDiff_k(\hat{P}_1,\hat{P})$.
\end{enumerate}
In coordinates, we have $(\sum_\sigma f_\sigma D_\sigma)^*=(-1)^{\abs{\sigma}
}\sum_\sigma D_\sigma\circ f_\sigma$ for a scalar operator and $\norm
{\Delta_{ij}}^*=\norm{\Delta_{ji}^*}$ for a matrix one.
\end{example}

\begin{example}
The choice $Q=Q_0\otimes\CL{p}^p$, with $Q_0$ a zero-curvature
representation, leads to the $p$-th \emph{gauge complex} (cf.\
\cite{Marvan1}).
\end{example}

\begin{example}
The zero term $E_0^{p,q}$ of \css consists of the complexes
\[0\xra{}E_0^{p,0}\xra{}E_0^{p,1}\xra{}\dotsb\xra{}E_0^{p,n}\xra{}0,
\quad p=0,1,2,\dots,\]
which are the horizontal de~Rham complexes with coefficients in $\CL{p}^p$.
\end{example}

\begin{example}\label{exmp:kr}
Take the \cm
\[Q=\bigoplus_p\Dv(\CL{p}^p)=\bigoplus_p\Hom_{\F}(\CL{}^1,\CL{p}^p).\]
The horizontal de~Rham complex with coefficient in $Q$ can be written as
\[0\xra{}\Dv\xra{}\Dv(\Lambda^1)\xra{}\Dv(\Lambda^2)\xra{}\dotsb\]
It is simply a calculation to verify that this is the complex introduced by
Krasil\cprime\-shchik \cite{Kras1}. The cohomology of this complex is called
the \emph{$\mathcal{C}$-co\-ho\-mol\-ogy}. It contains important invariants
of differential equations, in particular, symmetries and deformations of the
Cartan structure ($=$ recursion operators) (see
\cite{Kras1,Kras2,KrasKers1,KrasKers3,KrasKers2,KrasKers4}).
\end{example}

The calculation of the cohomology of the complexes from the last two
examples is our main concern in this paper.

\subsection{Compatibility complex}\label{subsec:cc}
Recall that a complex of \cd operators $\dotsb\xra{}P_{i-1}\xra{\Delta_i}P_i
\xra{\Delta_{i+1}}P_{i+1}\xra{}\dotsb$ is called \emph{formally exact}, if
the complex
\[\dotsb \xra{} \J^{k_i+k_{i+1}+l}(P_{i-1})
\xra{\varphi_{\Delta_i}^{k_i+k_{i+1}+l}} \J^{k_{i+1}+l}(P_i)
\xra{\varphi_{\Delta_{i+1}}^{k_{i+1}+l}}\J^l(P_{i+1}) \xra{} \dotsb,\]
with $\ord\Delta_j\le k_j$, is exact for any $l$.

Consider a \cd operator $\Delta\colon P_0\to P_1$ and the corresponding
\cm $\R_{\Delta}=\ker\varphi_{\Delta}^\infty$ (cf.\ example
\ref{exmp:ker}). Suppose that there exists a formally exact complex
\begin{equation}\label{eq:fe}
P_0 \xra{\Delta} P_1 \xra{\Delta_1} P_2 \xra{\Delta_2} P_3
\xra{\Delta_3} \dotsb.
\end{equation}
Then the cohomology of this complex coincides with the horizontal
cohomology with coefficients in $\R_{\Delta}$:
\begin{theorem}
\[\hH^i(\R_{\Delta})=H^i(P_{\bullet}).\]
\end{theorem}
\begin{proof}
Consider the following commutative diagram
\[\begin{CD}
@. \vdots @. \vdots @. \vdots @. \\
@.  @AAA  @AAA  @AAA  @. \\
0 @>>> \hL^2\otimes\Ji(P_0) @>>> \hL^2\otimes\Ji(P_1)
@>>> \hL^2\otimes\Ji(P_2) @>>> \dotsb \\
@.  @AA\hd A  @AA\hd A  @AA\hd A  @. \\
0 @>>> \hL^1\otimes\Ji(P_0) @>>> \hL^1\otimes\Ji(P_1)
@>>> \hL^1\otimes\Ji(P_2) @>>> \dotsb \\
@.  @AA\hd A  @AA\hd A  @AA\hd A  @. \\
0 @>>> \Ji(P_0) @>>> \Ji(P_1) @>>> \Ji(P_2) @>>> \dotsb \\
@.  @AAA  @AAA  @AAA  @. \\
@.  0  @. 0 @. 0  @.
\end{CD}\]
The horizontal maps are induced by the operators $\Delta_i$. All the
sequences are exact except for the terms in the left column and the
bottom row.  Now the standard spectral sequence arguments completes the
proof.
\end{proof}

Let us multiply the previous diagram by a \cm $Q$. This yields
\begin{equation}\label{eq:rd}
\hH^i(\R_{\Delta}\hatotimes Q)=H^i(P_{\bullet}\otimes Q),
\end{equation}
where $\R_{\Delta}\hatotimes Q=\projlim\R_{\Delta}^l\otimes Q$, with
$\R_{\Delta}^l=\ker\varphi_{\Delta}^{k+l}$, $\ord\Delta\le k$.

We can dualize our discussion. Namely, consider the diagram
\[\minCDarrowwidth=23pt\begin{CD}
@. \vdots @. \vdots @. \vdots @. \\
@.  @VVV  @VVV  @VVV  @. \\
0 @<<< \CDiff(P_0,\hL^{n-2}) @<<< \CDiff(P_1,\hL^{n-2})
@<<< \CDiff(P_2,\hL^{n-2}) @<<< \dotsb \\
@.  @VV\hd V  @VV\hd V  @VV\hd V  @. \\
0 @<<< \CDiff(P_0,\hL^{n-1}) @<<< \CDiff(P_1,\hL^{n-1})
@<<< \CDiff(P_2,\hL^{n-1}) @<<< \dotsb \\
@.  @VV\hd V  @VV\hd V  @VV\hd V  @. \\
0 @<<< \CDiff(P_0,\hL^n) @<<< \CDiff(P_1,\hL^n) @<<< \CDiff(P_2,\hL^n)
@<<< \dotsb \\
@.  @VVV  @VVV  @VVV  @. \\
@.  0  @. 0 @. 0  @.
\end{CD}\]

As above, we readily obtain
\[\hH^i(\R_{\Delta}^*)=H_{n-i}(\hat{P}_{\bullet})\]
and, more generally,
\begin{equation}\label{eq:rdst}
\hH^i(\R_{\Delta}^*\otimes Q)=H_{n-i}(\hat{P}_{\bullet}\otimes Q),
\end{equation}
where $\R_{\Delta}^*=\Hom(\R_{\Delta},\F)$. The homology in the right-hand
side of these formulae is the homology of the complex
\[\hat{P}_0 \xla{\Delta^*} \hat{P}_1 \xla{\Delta_1^*} \hat{P}_2
\xla{\Delta_2^*} \hat{P}_3 \xla{\Delta_3^*} \dotsb,\]
dual to the complex \eqref{eq:fe}.

Now, suppose we are given a \cd operator $\Delta\in\CDiff_k(P_0,P_1)$. How do
we find a formally exact complex of the form \eqref{eq:fe}? To this end
consider for each positive integer $k_1$ the mapping $\varphi_\Delta^{k+k_1}
\colon\J^{k+k_1}(P_0)\to\J^{k_1}(P_1)$. Without loss of generality it can be
assumed that for $k_1=0$ this mapping is surjective. Fix the integer $k_1>0$.
By definition, put $P_2=\coker\varphi_{\Delta}^{k+k_1}$ and let $\varphi_{
\Delta_1}$ be the natural projection $\J^{k_1}(P_1)\to P_2$. It is clear that
$\Delta_1\circ\Delta=0$. Further, starting from an integer $k_2$ and the
operator $\Delta_1$, we construct as above an operator $\Delta_2$, such that
$\Delta_2\circ\Delta_1=0$. Continuing this process, we obtain a complex of
the form \eqref{eq:fe}.

Is this complex formally exact? To settle the question, let us discuss the
notion of involutiveness of a \cd operator. Consider the Spencer
$\delta$-complex of the module $P_0$ (see Diagram \ref{eq:dcs})
\begin{equation}\label{eq:dsp}
0 \xra{}\hS^r\otimes P_0\xra{\bar{\delta}}\hL^1\otimes\hS^{r-1}\otimes P_0
\xra{\bar{\delta}}\hL^2\otimes\hS^{r-2}\otimes P_0\xra{\bar{\delta}}\dotsb
\end{equation}
Let $g^{k+l}\subset\hS^{k+l}\otimes P_0$ be the symbolic module of the
operator $\Delta$, i.e., $g^{k+l}=\ker\sigma(\Delta)$. It is easily
shown that the subcomplex of complex \eqref{eq:dsp}
\begin{equation}\label{eq:dspd}
0 \xra{} g^{k+l} \xra{\bar{\delta}} \hL^1\otimes g^{k+l-1}
\xra{\bar{\delta}} \hL^2\otimes g^{k+l-2} \xra{\bar{\delta}} \dotsb
\end{equation}
is well defined. Cohomology of this complex in the term $\hL^i\otimes g^{k+l
-i}$ is denoted by $\hH^{k+l,i}(\Delta)$ and is said to be \emph{horizontal
Spencer $\delta$-cohomology of the operator $\Delta$}. Note that $\hH^{k+l,0}
(\Delta)=\hH^{k+l,1}(\Delta)=0$. The operator $\Delta$ is called
\emph{involutive} (in the sense of Cartan), if $\hH^{k+l,i}(\Delta)=0$
for all $i\ge 0$.

\begin{theorem}
If the operator $\Delta$ is involutive, then the complex of the form
\eqref{eq:fe} constructed as described above is formally exact for all
positive integers $k_1$, $k_2$, $k_3$, \dots.
\end{theorem}
\begin{proof}
We must prove that the sequences
\[\hS^{k_{i-1}+k_i+l}\otimes P_{i-1} \xra{} \hS^{k_i+l}\otimes P_i \xra{}
\hS^l\otimes P_{i+1}\]
are exact for all $l\ge 1$. The proof is by induction on $i$ and $l$, with
the inductive step involving the standard spectral sequence arguments applied
to the commutative diagram
\[\minCDarrowwidth=24pt\begin{CD}
@.  0  @.  0  @.  0  @.  \\
@. @AAA  @AAA  @AAA @. \\
0 @>>> \hS^l\otimes P_{i+1} @>\bar{\delta}>>
\hL^1\otimes\hS^{l-1}\otimes P_{i+1} @>\bar{\delta}>>
\hL^2\otimes\hS^{l-2}\otimes P_{i+1} @>\bar{\delta}>> \dotsb \\
@. @AAA  @AAA  @AAA @. \\
0 @>>> \hS^{k_i+l}\otimes P_i @>\bar{\delta}>>
\hL^1\otimes\hS^{k_i+l-1}\otimes P_i @>\bar{\delta}>>
\hL^2\otimes\hS^{k_i+l-2}\otimes P_i @>\bar{\delta}>> \dotsb \\
@. @AAA  @AAA  @AAA @. \\
@. \vdots @. \vdots @. \vdots  @.  \\
@. @AAA  @AAA  @AAA @. \\
0 @>>> \hS^{K+l}\otimes P_0 @>\bar{\delta}>>
\hL^1\otimes\hS^{K+l-1}\otimes P_0 @>\bar{\delta}>>
\hL^2\otimes\hS^{K+l-2}\otimes P_0 @>\bar{\delta}>> \dotsb \\
@. @AAA  @AAA  @AAA @. \\
0 @>>> g^{K+l} @>\bar{\delta}>>
\hL^1\otimes g^{K+l-1} @>\bar{\delta}>>
\hL^2\otimes g^{K+l-2} @>\bar{\delta}>> \dotsb \\
@. @AAA  @AAA  @AAA @. \\
@.  0  @.  0  @.  0  @.
\end{CD}\]
where $K=k+k_1+k_2+\dots+k_i$.
\end{proof}

\begin{definition}
A formally exact complex of the form \eqref{eq:fe} is called the
\emph{compatibility complex} for the operator $\Delta$.
\end{definition}

For a discussion of the compatibility complex see \cite{BCGGG,Tsuji1,Verb5}.

The condition of involutiveness is not necessary for the existence of the
compatibility complex. Really, the $\delta$-Poincar\'e lemma (see, for
example, \cite{KLV,BCGGG}) says that for any \cd operator $\Delta$ there
exists an integer $l_0=l_0(m,n,k)$, where $m=\rank P$, such that
$\hH^{k+l,i}(\Delta)=0$ for $l\ge l_0$ and $i\ge 0$. Taking into
account this fact, we see from the proof of the previous theorem that
for sufficiently large integer $k_1$ the compatibility complex exists for
any \cd operator $\Delta$.

The following proposition is obvious.

\begin{proposition}\label{prop:comop}
In the compatibility complex each operator $\Delta_{i+1}$ is a compatibility
operator for the preceding operator $\Delta_i$, i.e., for any operator
$\nabla$ and an integer $l\ge 0$ such that $\nabla\circ\Delta_i=0$ and $\ord
\nabla\ge k_{i+1}-l$, there exists an operator $\square$ such that $\nabla^
{(l)}=\square\circ\Delta_{i+1}$, where $\nabla^{(l)}=\hj_l\circ\nabla$
is the $l$-th prolongation of $\nabla$.
\end{proposition}

\begin{example}
The de~Rham complex is the compatibility complex for the operator
$\hd\colon \F\to\hL^1$. The proof is trivial.
\end{example}

\begin{example}\label{exmp:dsd}
Fix an $\F$-linear scalar product of index $i$ on the module $\hL^1$.
For an integer $p\ge 1$ consider the operator
$\Delta=\hd{*}\hd\colon\hL^p\to\hL^ {n-p}$, where
$*\colon\hL^k\to\hL^{n-k}$ is the Hodge star operator. Let us show that
the complex
\[\hL^p \xra{\Delta} \hL^{n-p} \xra{\hd} \hL^{n-p+1} \xra{\hd}
\hL^{n-p+2} \xra{\hd} \dotsb \xra{\hd} \hL^n \xra{} 0\]
is the compatibility complex for the operator $\Delta$. Indeed, we must
prove that the image of the map
$\sigma(\Delta)\colon\hS^{l+2}\otimes\hL^p \to\hS^l\otimes\hL^{n-p}$
coincides with the image of the map $\sigma(\hd)
\colon\hS^{l+1}\otimes\hL^{n-p-1}\to\hS^l\otimes\hL^{n-p}$ for all
$l\ge 0$.  Since $\Delta*=\hd*\hd*=\hd(*\hd*+(-1)^{pn+n+i}\hd)$, it is
sufficient to show that the map
$\sigma(*\hd*+(-1)^{pn+n+i}\hd)\colon\hS^{l+1}\otimes\hL^{n-p-1}\to
\hS^l\otimes\hL^{n-p}$ is an epimorphism. Take an element
$\xi\in\hL^1$. One has
$\sigma(\hd)(\xi^{l+1}\otimes\omega)=l\xi^l\otimes\xi\wedge\omega$.
Denote by $A_{n-p-1}\colon\hL^{n-p-1}\to\hL^{n-p}$ the map of exterior
multiplication by $\xi$. We have $A_{n-p}^*=(-1)^{pn+n+i}*A_{p-1}*$,
so it will suffice to check that the map $A_{n-p-1}+A_{n-p}^*$ is an
epimorphism.  But this is obvious: $\hL^{n-p}=\im A_{n-p-1}\oplus(\im
A_{n-p-1})^{\bot}=\im A_{n-p-1}\oplus(\ker A_{n-p})^{\bot} =\im
A_{n-p-1}\oplus\im A_{n-p}^*=\im(A_ {n-p-1}+A_{n-p}^*)$\footnote{I
thank D. Gessler for drowning my attention to this example (cf.\
\cite{Gessler2}).}.
\end{example}

\section{Applications to computing cohomological invariants \\
         of systems of differential equations}\label{sec:appl}

\subsection{Main theorems}\label{subsec:main}
Let $\E=\{F=0\}$, $F\in P_1$, be an equation,
\[P_0=\varkappa \xra{\ell_F} P_1 \xra{\Delta_1} P_2 \xra{\Delta_2} P_3
\xra{\Delta_3} P_4 \xra{\Delta_4} \dotsb\]
the compatibility complex for the operator of universal linearization, and
\[\hat{P_0}=\hat{\varkappa}\xla{\ell_F^*}\hat{P_1}\xla{\Delta_1^*}\hat{P}_2
\xla{\Delta_2^*}\hat{P}_3\xla{\Delta_3^*}\hat{P}_4\xla{\Delta_4^*}\dotsb\]
the dual complex. Take a \cm $Q$.
\begin{theorem}
$\hH^i(\Dv(Q))=H^i(P_{\bullet}\otimes Q), \ \
\hH^{n-i}(\CL{}^1\otimes Q)=H_i(\hat{P}_{\bullet}\otimes Q).$
\end{theorem}
\begin{proof}
The statement follows immediately from \eqref{eq:rd}, \eqref{eq:rdst} and
Proposition \ref{prop:ker}.
\end{proof}

Let $Q=\CL{p}^p$. The previous theorem gives, first, a method for computing
of the cohomology groups $\hH^i(\Dv(\CL{p}^p))$, which are \kcc groups (see
Example \ref{exmp:kr}). Second, since the term $E_1^{p,q}=\hH^q(\CL{p}^p)$ of
\css is a direct summand in the cohomology group $\hH^q(\CL{}^1\otimes\CL
{p-1}^{p-1})$, we have a description for the first term of \css. Thus:

\begin{corollary}
$\hH^i(\Dv(\CL{p}^p))=H^i(P_{\bullet}\otimes\CL{p}^p)$.
\end{corollary}
\begin{corollary}
\label{cor:cssc}
The term $E_1^{p,q}$ of \css is the skew-symmetric part of the group
$H_{n-q}(\hat{P}_{\bullet}\otimes\CL{p-1}^{p-1})$.
\end{corollary}

It is useful to describe the isomorphisms given by these corollaries in an
explicit form. We discuss \css, the case of \kcc is similar.

Consider an operator $\nabla\in\CDiff(\varkappa,\hL^q\otimes\CL{p-1}^{p-1})$
that represents an element of $E_1^{p,q}$. This means that
\[\hd\circ\nabla=\nabla_1\circ\ell_F\]
for an operator $\nabla_1\in\CDiff(P_1,\hL^{q+1}\otimes\CL{p-1}^{p-1})$.
Applying the operator $\hd$ to both sides of this formula and using
Proposition \ref{prop:comop}, we get
\[\hd\circ\nabla_1=\nabla_2\circ\Delta_1\]
for an operator $\nabla_2\in\CDiff(P_2,\hL^{q+2}\otimes\CL{p-1}^{p-1})$.
Continuing this process, we obtain operators $\nabla_i\in\CDiff(P_i,\hL^{q+i}
\otimes\CL{p-1}^{p-1})$, $i=1,2,\dots,n-q$, such that
\[\hd\circ\nabla_{i-1}=\nabla_i\circ\Delta_{i-1}.\]
For $i=n-q$ this formula means that the operator $\nabla_{n-q}\in\CDiff(P_{n
-q},\hL^n\otimes\CL{p-1}^{p-1})$ represents an element of the module $\hat{P}
_{n-q}\otimes\CL{p-1}^{p-1}$ that lies in the kernel of the operator $\Delta_
{n-q-1}^*$. This is the element that gives rise to the homology class in
$H_{n-q}(\hat{P}_{\bullet}\otimes\CL{p-1}^{p-1})$ corresponding to the
chosen element of $E_1^{p,q}$.

It follows from our results that if there is an integer $k$ such that
$P_k=P_{k+1}=P_{k+2}=\dots=0$, i.e., the compatibility complex has the
form
\[P_0=\varkappa \xra{\ell_F} P_1 \xra{\Delta_1} P_2 \xra{\Delta_2} P_3
\xra{\Delta_3} \dotsb \xra{\Delta_{k-2}} P_{k-1} \xra{} 0,\]
then
\begin{enumerate}
\item $E_1^{p,q}=0$\quad for $p>0$ and $q\le n-k$,
\item $H^i(\Dv(\CL{p}^p))=0$\quad for $i\ge k$.
\end{enumerate}
This result is known as the \emph{$k$-line theorem}.

What are the values of the integer $k$ for differential equations encountered
in mathematical physics? The existence of a compatibility operator $\Delta_1$
is usually due to the existence of dependencies between the equations under
consideration: $\Delta_1(F)=0$. The majority of systems that occur in
practice consist of independent equations and for them $k=2$. In this case
the \emph{two-line theorem} holds:

\begin{theorem}[the two-line theorem]\label{thm:twoline}
Let an equation $\E$ be such that the compatibility complex for
$\ell_F$ is of length two. Then:
\begin{enumerate}
\item $E_1^{p,q}=0$\quad for $p>0$ and $q\le n-2$,
\item $E_1^{p,n-1}\subset\ker(\ell_F^*)_{\CL{p-1}^{p-1}}$\quad for $p>0$,
\item $E_1^{p,n}\subset\coker(\ell_F^*)_{\CL{p-1}^{p-1}}$\quad for $p>0$,
\item $H^i(\Dv(\CL{p}^p))=0$\quad for $i\ge 2$,
\item $H^0(\Dv(\CL{p}^p))=\ker(\ell_F)_{\CL{p}^p}$,
\item $H^1(\Dv(\CL{p}^p))=\coker(\ell_F)_{\CL{p}^p}$. \qed
\end{enumerate}
\end{theorem}
Further, we meet with the case $k>2$ in gauge theories, when the
dependencies $\Delta_1(F)=0$ are given by the second Noether theorem.
For usual irreducible gauge theories, like electromagnetism,
Yang\,-\,Mills models, and Einstein's gravity, the Noether identities
are independent, so that the operator $\Delta_2$ is trivial and, thus,
$k=3$. Finally, for an $L$-th stage reducible gauge theory, one has
$k=L+3$.

\subsection{Example: Evolution equations}\label{subsec:exeveq}
Consider an evolution equation $\E=\{F=u_t-f(x,t,u_i)=0\}$, with
independent variables $x,t$ and dependent variable $u$; $u_i$ denotes
the set of variables corresponding to derivatives of $u$ with respect
to $x$.

Natural coordinates for $\Ei$ are $(x,t,u_i)$. The total derivatives
operators $D_x$ and $D_t$ on $\Ei$ have the form
\[D_x=\dd{}{x}+\sum_i u_{i+1}\dd{}{u_i},\quad
D_t=\dd{}{t}+\sum_i D_x^i(f)\dd{}{u_i}.\]
The operator of universal linearization is given by
\[\ell_F=D_t-\ell_f=D_t-\sum_i\dd{f}{u_i}D_x^i.\]
The adjoint of $\ell_F$ is
\[\ell_F^*=-D_t-\ell_f^*=D_t-\sum_i (-1)^i D_x^i\circ\dd{f}{u_i}.\]

Clearly, for an evolution equation the two-line theorem holds, hence
\kcc $\hH^q(\Dv(\CL{p}^p))$ is trivial for $q\ge 2$ and the
first term $E_1^ {p,q}$ of \css is trivial for $q\ne 1,2$ and $p>0$.
Now, assume that the order of the equation $\E$ is greater than or
equal to $2$, i.e., $\ord\ell_f \ge 2$. Then one has more:

\begin{theorem} For any evolution equation of order $\ge2$, one has
\begin{enumerate}
\item $\hH^0(\Dv(\CL{p}^p))=0$\quad for $p\ge 2$,
\item $E_1^{p,1}=0$\quad for $p\ge 3$.
\end{enumerate}
\end{theorem}

\begin{proof}
It follows from Theorem \ref{thm:twoline} that
$\hH^0(\Dv(\CL{p}^p))=\ker( \ell_F)_{\CL{p}^p}$ and
$E_1^{p,1}=\ker(\ell_F^*)_{\CL{p-1}^{p-1}}$. Hence to prove the theorem
it suffices to check that equations
\begin{equation}\label{eq:do}
(D_t-\ell_f)(\omega)=0\quad\text{and}\quad(D_t+\ell_f^*)(\omega)=0,
\end{equation}
with $\omega\in\CL{p}^{p}$, has no nontrivial solutions for $p\ge 2$.

To this end, consider the symbol of \eqref{eq:do}. Denote $\sigma(D_x)=
\theta$. The symbol of $\ell_F$ has the form
$\sigma(\ell_f)=f\theta^k$, $k \ge 2$. An element $\omega\in\CL{p}^{p}$
can be identified with a multilinear \cd operator, so the symbol of
$\omega$ is a homogeneous polynomial in $p$ variables
$\sigma(\omega)=\delta(\theta_1,\dots,\theta_p)$.  Either of the two
equations \eqref{eq:do} yields
\[[f(\theta_1^k+\dots+\theta_p^k)\pm
f(\theta_1+\dots+\theta_p)^k](\delta)=0.\]
The conditions $k\ge 2$ and $p\ge 2$ obviously imply that $\delta=0$. This
completes the proof.
\end{proof}

\begin{remark}
In the work \cite{Gessler1} this proof has been generalized for
determined systems of evolution equations with arbitrary number of
independent variables.
\end{remark}

\subsection{Example: Abelian $p$-form theories}\label{subsec:exapf}

Let $M$ be a (pseudo-)Riemannian manifold and $\pi\colon E\to M$ the $p$-th
exterior power of the cotangent bundle over $M$, so that a section of $\pi$
is a $p$-form on $M$. Evidently, on the jet space $J^\infty(\pi)$ there
exists a unique horizontal form $A\in\hL^p(J^\infty(\pi))$ such that $j_
\infty^*(\omega)(A)=\omega$ for all $\omega\in\Lambda^p(M)$. Consider the
equation $\E=\{F=0\}$, with $F=\hd{*}\hd A$. Our aim is to calculate the
terms of \css $E_1^{i,q}$ for $q\le n-2$. We shall assume that $1\le p<n-1$
and that the manifold $M$ is topologically trivial.

Obviously, we have $P_0=\varkappa=\hL^p$, $P_1=\hL^{n-p}$, and $\ell_F=\hd{*}
\hd\colon\hL^p\to\hL^{n-p}$. Taking into account Example \ref{exmp:dsd}, we
see that the compatibility complex for $\ell_F$ has the form
\begin{equation}\label{eq:ccpf}
\begin{CD}
\hL^p @>\ell_F>> \hL^{n-p} @>\hd>> \hL^{n-p+1} @>\hd>> \dotsb
@>\hd>> \hL^n @>>> 0 \\
@|   @|   @|   @.  @|  \\
P_0 @. P_1 @. P_2 @. @. P_{k-1}
\end{CD}
\end{equation}
Thus $k=p+2$ and the $k$-line theorem yields $E_1^{i,q}=0$ for $i>0$ and $q<
n-p-1$. Since \css converges to the de~Rham cohomology of $\Ei$, which is
trivial, we also get $E_1^{0,q}=0$ for $0<q<n-p-1$, and $\dim E_1^{0,0}=1$,
i.e., $\hH^1=\hH^2=\dots=\hH^{n-p-2}=0$ and $\dim\hH^0=1$. Next, consider the
terms $E_1^{i,q}$ for $n-p-1\le q<2(n-p-1)$ and $i>0$. In view of Corollary
\ref{cor:cssc} one has
\[E_1^{i,q}\subset\hH^{q-(n-p-1)}(\CL{i-1}^{i-1})=E_1^{i-1,q-(n-p-1)}\]
because the complex dual to the compatibility complex \eqref{eq:ccpf} has the
form
\[\begin{CD}
\hL^{n-p} @<\ell_F^*<< \hL^p @<\hd<< \hL^{p-1} @<\hd<< \dotsb
@<\hd<< \F @<<< 0. \\
@|   @|   @|   @.  @|  \\
\hat{P}_0 @. \hat{P}_1 @. \hat{P}_2 @. @. \hat{P}_{p+1}
\end{CD}\]
(Throughout, it is assumed that $q\le n-2$.) Thus we obtain $E_1^{i,q}=0$ for
$n-p-1<q<2(n-p-1)$, $i>0$ and $\dim E_1^{1,n-p-1}=1$. Again, taking into
account that the spectral sequence converges to the trivial cohomology, we
get $E_1^{0,q}=0$ for $n-p-1<q<2(n-p-1)$ and $\dim E_1^{0,n-p-1}=1$. In
addition, the map $d_1^{0,n-p-1}\colon E_1^{0,n-p-1}\to E_1^{1,n-p-1}$
is an isomorphism. Explicitly, one readily obtains that the
one-dimensional space $E_1^{0,n-p-1}$ is generated by the element $*\hd
A\in\hL^{n-p-1}$ and the map $d_1^{0,n-p-1}$ takes this element to the
operator $*\hd\colon\varkappa= \hL^p\to\hL^{n-p-1}$, which generates
the space $E_1^{1,n-p-1}$.

Further, let us consider the terms $E_1^{i,q}$ for $2(n-p-1)\le q<3(n-p-1)$.
Arguing as before, we see that all these terms vanish unless $q=2(n-p-1)$ and
$i=0,1,2$, with $\dim E_1^{1,2(n-p-1)}=1$ and $\dim E_1^{i,2(n-p-1)}\le 1$,
$i=0,2$. To compute the terms $E_1^{i,2(n-p-1)}$ for $i=0$ and $i=2$, we have
to consider two cases: $n-p-1$ is even and $n-p-1$ is odd (see Diagram
\ref{pic}).

\begin{figure}
\begin{picture}(340,200)
\put(70,2){$n-p-1$ is even}
\put(250,2){$n-p-1$ is odd}
\put(60,15){\vector(1,0){100}}
\put(60,15){\vector(0,1){185}}
\put(240,15){\vector(1,0){100}}
\put(240,15){\vector(0,1){185}}
\put(53,195){$q$}
\put(156,5){$i$}
\put(233,195){$q$}
\put(336,5){$i$}
\multiput(78,15)(18,0){5}{\line(0,1){185}}
\multiput(258,15)(18,0){5}{\line(0,1){185}}
\put(55,180){\llap{$3(n-p-1)$}}
\put(55,125){\llap{$2(n-p-1)$}}
\put(55,70){\llap{$n-p-1$}}
\put(235,180){\llap{$3(n-p-1)$}}
\put(235,125){\llap{$2(n-p-1)$}}
\put(235,70){\llap{$n-p-1$}}
\put(69,19){\circle*{3}}
\multiput(69,72)(0,55){3}{\circle*{3}}
\multiput(87,72)(0,55){3}{\circle*{3}}
\multiput(72,72)(0,55){3}{\vector(1,0){12}}
\put(249,19){\circle*{3}}
\multiput(249,72)(18,55){3}{\circle*{3}}
\multiput(267,72)(18,55){3}{\circle*{3}}
\multiput(252,72)(18,55){3}{\vector(1,0){12}}
\end{picture}
\caption{}
\label{pic}
\end{figure}

In the first case, the map $d_1^{1,2(n-p-1)}\colon E_1^{1,2(n-p-1)}\to
E_1^ {2,2(n-p-1)}$ is trivial. Indeed, the operator $(*\hd
A)\wedge*\hd\colon \varkappa=\hL^p\to\hL^{2(n-p-1)}$, which generates
the space $E_1^{1,2(n-p- 1)}$, under the mapping $d_1^{1,2(n-p-1)}$ is
the antisymmetrization of the operator
$\omega_1\times\omega_2\mapsto(*\hd\omega_1)\wedge(*\hd\omega_2)$,
$\omega_i\in\varkappa=\hL^p$. But this operator is symmetric, so that
$d_1^ {1,2(n-p-1)}=0$. Consequently, $E_1^{2,2(n-p-1)}=0$ and $\dim
E_1^{0,2(n-p- 1)}=1$. This settles the case when $n-p-1$ is even.

In the case when $n-p-1$ is odd, the operator $\omega_1\times\omega_2\mapsto
(*\hd\omega_1)\wedge(*\hd\omega_2)$ is skew-symmetric, hence the map
$d_1^{1, 2(n-p-1)}$ is an isomorphism. Thus, $\dim E_1^{2,2(n-p-1)}=1$
and $E_1^{0, 2(n-p-1)}=0$.

Continuing this line of reasoning, we obtain the following result.

\begin{theorem} For $i=q=0$ one has $\dim E_1^{0,0}=1$. If either or
both $i$ or $q$ are positive, there are two cases:

\begin{enumerate}
\item if $n-p-1$ is even then $\dim E_1^{i,q}=
\begin{cases}
1&\text{for $i=l(n-p-1)$ and $q=0$, $1$}, \\
0&\text{otherwise};
\end{cases}$
\item if $n-p-1$ is odd then $\dim E_1^{i,q}=
\begin{cases}
1&\text{for $i=l(n-p-1)$ and $q=l-1$, $l$}, \\
0&\text{otherwise};
\end{cases}$
\end{enumerate}
here $1\le l<\dfrac{n-1}{n-p-1}$.
\end{theorem}

In other words, let $\bar{\mathcal{A}}$ be the exterior algebra generated by
two forms: $\omega_1=*\hd A\in\hL^{n-p-1}$ and $\omega_2=\hd_1(\omega_1)=*\hd
\in\hL^{n-p-1}\otimes\CL{}^1$; then we see that the space $\bigoplus_{i,q\le
n-2}E_1^{i,q}$ is isomorphic to the subspace of $\bar{\mathcal{A}}$
containing no forms of degree $q>n-2$. In particular, for horizontal
cohomology with trivial coefficients this result agrees with that of
\cite{HennKnaSchom} obtained previously by means of the Koszul\,-\,Tate
resolution.

\begin{remark}
One can consider the quotient equation under the action of gauge symmetries.
The techniques described here allow one to compute the terms $E_1^{i,q}$,
$q<n-1$, for the quotient equation. This will be considered elsewhere.
\end{remark}

\end{document}